\documentclass[12pt]{amsart}
\usepackage{amssymb}

\newtheorem{thm}{Theorem}[section]
\newtheorem{Prop}[thm]{Proposition}
\newtheorem{lemma}[thm]{Lemma}

\newtheorem{co}[thm]{Corollary}

\newcommand{\I}{{\text{Im}}}

\newcommand{\bc}{{\mathbb C}}

\newcommand{\bh}{{\mathbb H}}

\newcommand{\br}{{\mathbb R}}

\newcommand{\ra}{\rightarrow}

\newcommand{\hh}{H_{\bh}^}

\newenvironment{pf}{\begin{trivlist}\item[]{\bf Proof:\ }}
{\mbox{}\hfill\rule{.08in}{.08in}\end{trivlist}}

\begin{document}
\title[Rigidity of quaternionic  space]
{Local Rigidity  in quaternionic hyperbolic space}
\author{Inkang Kim and Pierre Pansu}
\date{}
\maketitle

\begin{abstract}
In this note, we study deformations of quaternionic hyperbolic
lattices in larger quaternionic hyperbolic spaces and prove local
rigidity results. On the other hand, surface groups are shown to be
more flexible in quaternionic hyperbolic plane than in complex
hyperbolic plane.
\end{abstract}
\footnotetext[1]{I. Kim gratefully acknowledges the partial support
of KRF Grant (0409-20060066) and a warm support of IHES during his
stay.} \footnotetext[2]{P. Pansu, Univ Paris-Sud, Laboratoire de
Math\'ematiques d'Orsay, Orsay, F-91405} \footnotetext[3]{\hskip
47pt CNRS, Orsay, F-91405.}

\section{Introduction}

\subsection{4-dimensional lattices}

Lattices in $Sp(n,1)$, $n\geq 2$, when map\-ped to $Sp(m,1)$, cannot
be deformed. This follows from K. Corlette's archimedean
superrigidity theorem, \cite{C2}. What about lattices in $Sp(1,1)$,
i.e. in 4-dimensional hyberbolic space ?

In this note we prove local rigidity of uniform lattices of $Sp(1,1)$ when mapped to $Sp(2,1)$. In complex hyperbolic geometry, such rigidity results were first discovered by D. Toledo, \cite{To}. In \cite{G1,GM}, W. Goldman and J. Millson gave a local explanation of this phenomenon. Our main result is an exact quaternionic analogue of theirs.

Start with a uniform lattice $\Gamma$ in $Sp(1,1)$. There is an easy
manner to deform the embedding $\rho_0 :\Gamma\ra Sp(1,1)\ra
Sp(2,1)$. Indeed, since $Sp(2,1)$ contains $Sp(1,1)\times Sp(1)$, it
also contains many copies of $Sp(1,1)\times U(1)$. If $H^1
(\Gamma,\mathbb{R})\not=0$, which happens sometimes (see
\cite{Mill}), the trivial representation $\Gamma\to U(1)$ can be
continuously deformed to a nontrivial representation $\rho_1$. All
such representations give rise to actions on quaternionic hyperbolic
plane which stabilize a quaternionic line. Therefore, only
deformations which do not stabilize any quaternionic line should be
of interest.

\begin{thm}
\label{4} Let $\Gamma \subset Sp(1,1)$ be a lattice. Embed $\Gamma$
into $Sp(2,1)$ as a subgroup which stabilizes a quaternionic line.

If $\Gamma$ is uniform in $Sp(1,1)$, then every small deformation of
$\Gamma$ in $Sp(2,1)$ again stabilizes a quaternionic line.

If $\Gamma$ is non uniform in $Sp(1,1)$, then every small
deformation of $\Gamma$ in $Sp(2,1)$ preserving parabolics again
stabilizes a quaternionic line.
\end{thm}

Toledo's theorem inaugurated a series of global rigidity results by A. Domic, D. Toledo, \cite{DT}, K. Corlette, \cite{C0}, M. Burger, A. Iozzi and A. Wienhard, \cite{BIW}. By global rigidity, we mean the following : a certain characteristic number of representations, known as Toledo
invariant, is maximal if and only if the representation stabilizes a
totally geodesic complex hypersurface. It is highly expected that
such a global rigidity should hold in quaternionic hyperbolic
spaces, but we have been unable to prove it. Note that since $Sp(1,1)=Spin(4,1)^0$, there exist uniform lattices in $Sp(1,1)$ which are isomorphic to Zariski dense subgroups of $Sp(4,1)$, see section \ref{bend}.

\textbf{Question}. Let $\Gamma \subset Sp(1,1)$ be a uniform
lattice. Embed $\Gamma$ into $Sp(3,1)$. Can one deform $\Gamma$ to a
Zariski dense subgroup?

\subsection{3-dimensional lattices}

Uniform lattices in 3-dimensional real hyperbolic space can
sometimes be deformed nontrivially in 4-dimensio\-nal real
hyperbolic space, see \cite{Th}, chapter 6, or \cite{Apa}.
Nevertheless, when they act on quaternionic plane, all small
deformations stabilize a quaternionic line, although the action on
this line can be deformed non trivially.

\begin{thm}
\label{3} Let $\Gamma\subset Spin(3,1)^0$ be a lattice. Embed
$Spin(3,1)^0$ into $Spin(4,1)^0 =Sp(1,1)$ and then into $Sp(2,1)$ in
the obvious manner. This produces a discrete subgroup of $Sp(2,1)$
stabilizing a quaternionic line.

If $\Gamma$ is uniform in $Spin(3,1)^0$, then every small
deformation of $\Gamma$ in $Sp(2,1)$ again stabilizes a quaternionic
line.

If $\Gamma$ is non uniform in $Spin(3,1)^0$, then every small
deformation of $\Gamma$ preserving parabolics again stabilizes a
quaternionic line.
\end{thm}

If the assumption on parabolics is removed, nonuniform lattices in
$Spin(3,1)^0$ can be deformed within $Spin(3,1)^0$, see \cite{Th},
chapter 5.

\textbf{Question}. Let $\Gamma$ be a non uniform lattice in
$Spin(3,1)^0$. Map it to $Sp(2,1)$ via $Spin(4,1)^0 =Sp(1,1)$. Can
one deform $\Gamma$ to a Zariski-dense subgroup?

\subsection{2-dimensional lattices}

Uniform lattices in real hyperbolic pla\-ne, when mapped to
$SU(2,1)$ using the embedding $SO(2,1)\ra SU(2,1)$, can be deformed
to discrete Zariski-dense subgroups of $SU(2,1)$. On the other hand,
lattices mapped via $SU(1,1)$ and $SU(2,1)$ are more rigid, as shown
by W. Goldman and J. Millson, \cite{GM}. This fact has been recently
extended to higher rank groups by M. Burger, A. Iozzi and A.
Wienhard, \cite{BIW}.

It turns out that this form of rigidity of surface groups does not
apply to the group $Sp(2,1)$.

\begin{thm}
\label{2} Let $\Gamma$ be the fundamental group of a closed surface
of genus $>1$.  \begin{enumerate}
  \item View $\Gamma$ as a uniform lattice in $SO(2,1)$. Map $SO(2,1) \ra  Sp(2,1)$. This gives rise to a representation into $Sp(2,1)$ which can be deformed to a discrete Zariski-dense representation.
  \item View $\Gamma$ as a uniform lattice in $SU(1,1)$. Map $SU(1,1) \ra Sp(1,1)\ra Sp(2,1)$. This gives rise to a representation into $Sp(2,1)$ fixing a quaternionic line. Then there exists small deformations which do not stabilize any quaternionic line.
\end{enumerate}
\end{thm}

Whereas in the first case, explicit examples of deformations are
provided by Thurston's bending construction, the existence of
Zariski dense deformations in the second case follows from rather
general principles. It would be interesting to visualize some of
them.

\subsection{Plan of the paper}

Section \ref{rel} gives a cohomological criterion for non Zariski
dense sugroups to remain non Zariski dense after deformation. The
necessary cohomology vanishing is obtained in section \ref{rep}.
Theorem \ref{4} is proved in section \ref{4u}, Theorem \ref{3} in
section \ref{low}. The statements for nonuniform lattices are proved
in section \ref{nonuniform}. Section \ref{bend} describes how
lattices in Lie subgroups can sometimes be bent to become Zariski
dense. The proof of Theorem \ref{2} is completed in section
\ref{fuchsian}. We end with a remark on non Zariski dense discrete
subgroups in section \ref{discrete}.

\section{A relative Weil theorem}
\label{rel}

Let $\Gamma$ be a finitely generated group, and $G$ be a Lie group
with Lie algebra $\mathfrak{g}$. The {\em character variety}
$\chi(\Gamma,G)$ is the quotient of the space $Hom(\Gamma,G)$ of
homomorphisms of $\Gamma$ to $G$ by the action of $G$ by
postcomposing homomorphisms with inner automorphisms. In
\cite{Weil}, A. Weil shows that a sufficient condition for a
homomorphism $\rho:\Gamma\ra G$ to define an isolated point in the
character variety is that the first cohomology group $H^1
(\Gamma,\mathfrak{g}_\rho)$ vanishes. In this section, we state a
relative version of Weil's theorem.

Let $H\subset G$ be an algebraic subgroup of $G$. Let
$\chi(\Gamma,H,G)\subset\chi(\Gamma,G)$ be the set of conjugacy
classes of homomorphisms $\Gamma\to G$ which fall into conjugates of
$H$. In other words, $\chi(\Gamma,H,G)$ is the set of $G$-orbits of
elements of $Hom(\Gamma,H)\subset Hom(\Gamma,G)$. If $\rho\in
Hom(\Gamma,H)$, the representation $\mathfrak{g}_\rho =ad\circ\rho$
on the Lie algebra $\mathfrak{g}$ of $G$ leaves the Lie algebra
$\mathfrak{h}$ of $H$ invariant, and thus defines a quotient
representation, which we shall denote by
$\mathfrak{g}_{\rho}/\mathfrak{h}_{\rho}$.

\begin{Prop}
\label{relative} Let $H\subset G$ be real Lie groups, with Lie
algebras $\mathfrak{h}$ and $\mathfrak{g}$. Let $\Gamma$ be a
finitely generated group. Let $\rho:\Gamma\ra H$ be a homomorphism.
Assume that $H^1
(\Gamma,\mathfrak{g}_{\rho}/\mathfrak{h}_{\rho})=0$. Then
$\chi(\Gamma,H,G)$ is a neighborhood of the $G$-conjugacy class of
$\rho$ in $\chi(\Gamma,G)$. In other words, homomorphisms $\Gamma\ra
G$ which are sufficiently close to $\rho$ can be conjugated into
$H$.
\end{Prop}

\begin{pf}
$Hom(\Gamma,G)$ is topologized as a subset of the space $G^{\Gamma}$
of arbitrary maps $\Gamma\ra G$. Let $\Phi:G^{\Gamma}\ra
G^{\Gamma\times\Gamma}$ be the map which to a map $f:\Gamma\ra G$
associates $\Phi(f):\Gamma\times\Gamma\ra G$ defined by
\begin{eqnarray*}
\Phi(f)(\gamma,\gamma')=f(\gamma\gamma'^{-1})f(\gamma)f(\gamma').
\end{eqnarray*}
In other words, a map $f\in G^{\Gamma}$ is a homomorphism if and
only if $\Phi(f)=1$.

Consider the map $\Psi:G\times H^{\Gamma}\ra G^{\Gamma}$ which sends
$g\in G$ and $f:\Gamma\ra H$ to the map $\Psi(g,f):\Gamma\to G$
defined by
\begin{eqnarray*}
\Psi(g,f)(\gamma)=g^{-1}f(\gamma)g.
\end{eqnarray*}
We need prove that the image of $\Psi$ contains a neighborhood of
$\rho$ in $\Phi^{-1}(1)$.

The cohomological assumption gives information on the differentials
of $\Phi$ and $\Psi$. The differential $D_{\rho}\Phi$ is equal to
$-d_1$ where $d_1$ denotes the coboundary
$C^{1}(\Gamma,\mathfrak{g}_{\rho})\ra
C^{2}(\Gamma,\mathfrak{g}_{\rho})$. The differential of $\Psi$ at
$g=e$ and $f=\rho$ is given by
\begin{eqnarray*}
D_{(e,\rho)}\Psi(v,\eta)=-d_0 v +\eta,
\end{eqnarray*}
where $d_0$ denotes the coboundary
$C^{0}(\Gamma,\mathfrak{g}_{\rho})\ra
C^{1}(\Gamma,\mathfrak{g}_{\rho})$. Since, for all $f\in
H^{\Gamma}$,
$\Phi(\Psi(g,f))(\gamma,\gamma')=g^{-1}\Phi(f)(\gamma,\gamma')g$,
$D_{\rho}\Phi\circ D_{(e,\rho)}\Psi=0$. Conversely, if we assume
that $H^1 (\Gamma,\mathfrak{g}_{\rho}/\mathfrak{h}_{\rho})=0$, any
$\theta\in C^1 (\Gamma,\mathfrak{g}_{\rho})$ such that
$D_{\rho}\Phi(\theta)$ takes values in the subalgebra $\mathfrak{h}$
can be written $\theta=-d_0 v +\eta$ where $v\in \mathfrak{g}$ and
$\eta\in C^1(\Gamma,\mathfrak{h}_{\rho})$, i.e. $\theta$ belongs to
the image of $D_{(e,\rho)}\Psi$.

Clearly, $Hom(\Gamma,G)$ and $Hom(\Gamma,H)$ are real analytic
varieties. To analyze a neigborhood of $\rho$ in them, it is
sufficient to analyze real analytic of even formal curves
$t\mapsto\rho(t)$. In coordinates for $G$ (in which $H$ appears as a
linear subspace), such a curve admits a Taylor expansion
\begin{eqnarray*}
\rho(t)=\sum_{n=0}^{\infty}a_j t^j ,
\end{eqnarray*}
where $a_0 =\rho$ and for $j\geq 1$, $a_j \in
C^{1}(\Gamma,\mathfrak{g}_{\rho})$ is a $1$-cochain. Then
$\Phi(\rho(t))=1$ for all $t$. Expanding this as a Taylor series
gives
\begin{eqnarray*}
1&=&\Phi(\rho)+D_{\rho}\Phi(a_1)t
+(D_{\rho}\Phi(a_2)+D^{2}_{\rho}\Phi(a_1 ,a_1))t^2 +\cdots,
\end{eqnarray*}
which implies that
\begin{eqnarray*}
D_{\rho}\Phi(a_1)=0, \quad D_{\rho}\Phi(a_2)+D^{2}_{\rho}\Phi(a_1
,a_1)=0,\quad \ldots
\end{eqnarray*}
The first equation says that $a_1$ is a cocycle. So is $a_1$ mod
$\mathfrak{h}$, therefore there exist $v\in \mathfrak{g}$ and $b_1
\in Z^1 (\Gamma,\mathfrak{h}_\rho)$ such that $a_1 =-d_0 v+b_1$. Let
$t\mapsto g(t)$ be an analytic curve in $G$ with Taylor expansion
$g(t)=1+vt+\cdots$. Then the Taylor expansion of $\rho_1
(t)=g(t)^{-1}\rho(t)g(t)$ takes the form $\rho_1 (t)=1+b_1
t+\cdots$. In other words, up to conjugating, we arranged to bring
the first term of the expansion of $\rho(t)$ into $\mathfrak{h}$.

The second equation now reads
$D_{\rho}\Phi(a_2)+D^{2}_{\rho}\Phi(b_1 ,b_1)=0$. It implies that
$D_{\rho}\Phi(a_2)$ takes its values in $\mathfrak{h}$. Therefore
there exist $v'\in \mathfrak{g}$ and $b_2 \in Z^1
(\Gamma,\mathfrak{h}_{\rho})$ such that $a_2 =-d_0 v'+b_2$.
Conjugating $\rho_1 (t)$ by an analytic curve in $G$ with Taylor
expansion $1+v't^2 +\cdots$ kills $v'$ and replaces $a_2$ with $b_2$
in the expansion of $\rho_1 (t)$. Inductively, one can bring all
terms of the expansion of $\rho(t)$ into $\mathfrak{h}$. The
resulting curve belongs to $Hom(\Gamma,H)$. This shows that in a
neighborhood of $\rho$, $Hom(\Gamma,G)$ coincides with
$G^{-1}Hom(\Gamma,H)G$. Passing to the quotient, $\chi(\Gamma,H,G)$
coincides with $\chi(\Gamma,G)$ in a neighborhood of the conjugacy
class of $\rho$.
\end{pf}

\section{A cohomology vanishing result}
\label{rep}

\subsection{Preliminaries}
\label{pre}

For basic information on quaternionic hyperbolic space and surveys,
see \cite{KimJKMS,KP,Pa}.

We regard $\bh^n$ as a right module over $\bh$ by right
multiplication. Viewing $\bh=\bc \oplus j\bc= \bc^2$, left
multiplication by $\bh$ gives $\bc$-linear endomorphisms of $\bc^2$.
So $\bh^*=GL_1 \bh \subset GL_2\bc$. Similarly
$(x_1+iy_1+j(z_1+iw_1),\cdots,x_n+ i y_n + j(z_n+ i w_n))$ is
identified with $(x_1+iy_1,\cdots,x_n+i y_n; z_1+iw_1,\cdots,z_n+i
w_n)$ so that $\bh^n=\bc^{2n}$ and $GL_n\bh\subset GL_{2n}\bc$.

A $\bc$-linear map $\phi:\bh^n \ra \bh^n$ is $\bh$-linear exactly
when it commutes with $j:\phi(vj)=\phi(v)j$. Then it follows that if
$J=\left[
\begin{matrix}
     0 & - I_n \\
     I_n & 0 \end{matrix} \right],$
$$GL_n \bh=\{A\in GL_{2n}\bc: AJ=J\bar{A}\}.$$

Any element in $GL_n\bh$ can be written as $\alpha + j\beta$ where
$\alpha$ and $\beta$ are  ${2n}\times 2n$ complex matrices. If we
write a vector in $\bh^n$ in the form $X+jY$ where $X,Y\in \bc^n$,
the action of $\alpha + j\beta$ on it is
$$\alpha X -\overline \beta Y + j(\overline \alpha Y +\beta X).$$
So a matrix $\alpha+j\beta$ in $GL_n\bh$ corresponds to a matrix in
$GL_{2n}\bc$
\[ \left[\begin{matrix}
     \alpha  &  -\overline\beta \\
     \beta  &   \overline\alpha \end{matrix} \right].\]

In this paper, we fix a  quaternionic Hermitian form $\langle\
\rangle$ of signature $(n,1)$ on $\bh^{n+1}$ as $$\langle v, w
\rangle=\sum_{i=1}^n \bar{v_i}w_i- \overline{v_{n+1}}w_{n+1}.$$ Then
the Lie group $Sp(n,1,\bh)=Sp(n,1)$, which is the set of matrices
preserving this Hermitian form is
$$\{A\in GL_{n+1}\bh: A^*J'A=J'\},$$ where $J'=\left [
\begin{matrix}
  I_n & 0 \\
  0   & -1 \end{matrix} \right].$

It is easy to see that its Lie algebra $\mathfrak{sp}(n,1)$ is the
set of matrices of the form
$$\left[ \begin{matrix}
   \I\bh & Y \\
   X  &  \mathfrak{sp}(n-1,1) \end{matrix}\right],$$ where $Y +
   X^*J_n=0, X,Y\in \bh^n$ and
   $J_n=\left[ \begin{matrix}
         I_{n-1} & 0 \\
         0       & -1 \end{matrix}\right ].$
         So we get $$\mathfrak{sp}(n,1)=\I\bh \oplus \bh^n \oplus
\mathfrak{sp}(n-1,1).$$

Note that the adjoint action of the subgroup $\left [
\begin{matrix}
Sp(1) & 0\\
0  & Sp(n-1,1)\end{matrix}\right]$ preserves this decomposition. The
action on the $\bh^n$ component is the standard action,
$$Sp(n-1,1) \bh^n Sp(1)^{-1}.$$

Identifying $\bh^{n+1}$ with $\bc^{2n+2}$ as above, it is easy to
see that $Sp(n,1,\bh)$ is exactly equal to $U(2n,2)\cap
Sp(2n+2,\bc)$, i.e. to the set of unitary matrices satisfying
$AJ=J\bar{A}$. Indeed, the symplectic form with respect to the
standard basis of $\bc^{2n+2}$ is
\[ \left[
\begin{matrix}
                            0 &  A \\
                            -A & 0 \end{matrix} \right] \] and
                            $A=\left[ \begin{matrix}
                             I_n & 0 \\
                             0  & -1 \end{matrix} \right]$.

We will often complexify real Lie algebras. For any $M\in
\mathfrak{gl}(2n,\bc)$, one can write
$$M=\frac{1}{2}(M-J\bar{M}J)-i(\frac{1}{2}(iM+ i J\bar{M}J)).$$
So it is easy to see that
$\mathfrak{gl}(n,\bh)=\{A\in\mathfrak{gl}(2n,\bc):AJ=J\bar{A}\}$ is
complexified to $\mathfrak{gl}(2n,\bc)$. It is well-known that
$\mathfrak{u}(2n,2)\otimes_{\br} \bc=\mathfrak{gl}(2n+2,\bc)$ and
$\mathfrak{sp}(2n+2,\bc)\otimes_{\br}
\bc=\mathfrak{sp}(2n+2,\bc)\times \mathfrak{sp}(2n+2,\bc)$. From
these, we obtain that
$$\mathfrak{sp}(n,1)\otimes_{\br} \bc=\mathfrak{sp}(2n+2,\bc).$$
We are particulary interested in $$\mathfrak{sp}(1,1)\otimes_{\br}
\bc=\mathfrak{sp}(4,\bc).$$

The {\bf quaternionic hyperbolic} $n$-space $H^n_\bh$ in the unit
ball model is
$$\{(x_1,\cdots,x_n)|x_i \in \bh,\sum |x_i|^2 <1 \}.$$
It can be also described as a hyperboloid model
$$\{X\in \bh^{n+1}: \langle X, X \rangle=-1 \}/\sim$$ where
 $X \sim Y$ iff $X= YSp(1)$. Then the isometry group of $\hh n$ is
$PSp(n,1)$ which is a noncompact real semi-simple Lie group.

A point $X$ in the unit ball model can be mapped to $[X,1]$ in the
hyperboloid model. Then it is easy to see that the subgroup of the
form
$$\left[ \begin{matrix}
Sp(n-1) & 0 \\
0 & Sp(1,1) \end{matrix} \right]$$
stabilizes a quaternionic line $(0,0,\cdots,0,\bh)$ in the ball
model. In fact, we have
\begin{lemma}\label{stab}
The stabilizer of a quaternionic line $\{(0,\bh)\}$ in $Sp(2,1)$ is
of the form
$$\left [ \begin{matrix}
    Sp(1) &  0  \\
    0   &  Sp(1,1)
\end{matrix} \right].$$ Furthermore a parabolic element in $SO(4,1)=
Sp(1,1)$ stabilizing the quaternionic line $\{(0,\bh)\}$ is of the
form in $PSp(2,1)$
$$\left[ \begin{matrix}
  Sp(1)  &  0  \\
   0    & \left [ \begin{matrix}
            a  & \lambda-a \\
            a-\lambda  & 2\lambda-a  \end{matrix} \right]\end{matrix}\right]$$
where $a\geq 1$ is a positive real number, $\lambda \in Sp(1)$ with
${\text Re} \lambda= \frac{1}{a}$. These elements constitute the
parabolic elements in the center $\{(t,0)\}$ of the Heisenberg
group. A general parabolic element fixing a point $(0,1)$ at
infinity and not stabilizing the quaternionic line $\{(0,\bh)\}$, is
of the form
$$\left[ \begin{matrix}
     * & x & -x \\
     * & * & * \\
     * & * & * \end{matrix} \right],$$
with $x\neq 0$. These elements constitute the parabolic elements
which do not belong to the center of the Heisenberg group.
\end{lemma}
\begin{pf}
The quaternionic line $\{(0,\bh)\}$ in the hyperboloid model has
coordinate $(0,\bh,1)$. To fix this line, it is not difficult to see
that the matrix should have the form of $A=\left[ \begin{matrix}
  * & 0 & 0 \\
  * & * & 0 \\
  * & * & * \end{matrix} \right]$. Since its inverse $J'A^*J'$ also
fixes the quaternionic line, it should have the form as in the
claim.

Now to prove the second claim, note that the matrix should satisfy
the equation $\left[\begin{matrix}
                 Sp(1) &0 & 0 \\
                  0   & a & b \\
                  0   & c & d \end{matrix}\right]
(0,1,1)=\lambda (0,1,1)$ for $\lambda\in Sp(1)$. Also it should
satisfy $A^*J'A=J'$. From these, we obtain
$$a+b=\lambda$$
$$c+d=\lambda$$
$$|a|^2-|c|^2=|d|^2-|b|^2=1$$
$$\bar a b -\bar c d=0.$$
Then we get $\bar a(\lambda-a)-\bar c(\lambda-c)=0$. So $(\bar
a-\bar c)\lambda =|a|^2-|c|^2=1$, and we get $c=a-\lambda$. Now we
divide $A$ by $a$ since $a$ is nonzero. Note that $Aa^{-1}$
represents the same element in $PSp(2,1)$. Then we can assume that
$a$ is a positive real number, conjugating $A$ if necessary. The
fact that ${\text Re} \lambda=\frac{1}{a}$ follows from the other
two equations. So the result follows. In Heisenberg group
$\{(t,z)|t\in \I \bh, z\in \bh\}$, the center $\{(t,0)\}$ is the
(ideal) boundary of the quaternionic line $\{(0,\bh)\}$. So these
parabolic elements stabilizing the quaternionic line belong to the
center. See \cite{Kim}.

To prove the last claim, we just note that $A(0,1,1)=\lambda
(0,1,1)$ should be satisfied. The parabolic elements not stabilizing
the quaternionic line $\{(0,\bh)\}$ should have nonzero $x$ by the
first case.
\end{pf}

\subsection{Raghunathan's theorem}

In this section we collect information concerning finite dimensional
representations of $\mathfrak{so}(5,\bc)$, which will be necessary
for our main theorem. The basic theorem we will make use of is due
to M.S. Raghunathan, \cite{Ra}.

\begin{thm}
\label{raghunathan} Let $G$ be a connected semi-simple Lie group.
Let $\Gamma\subset G$ be a uniform irreducible lattice and
$\rho:(\Gamma \subset G)\ra Aut(E)$ a simple non-trivial linear
representation. Then $H^1(\Gamma;E)=0$ except possibly when
$\mathfrak g={\mathfrak{so}}(n+1,1)$ (resp. $\mathfrak
  g={\mathfrak{su}}(n,1))$ and the highest weight of $\rho$ is a
  multiple of the highest weight of the standard representation of
  ${\mathfrak {so}}(n+1,1)$ (resp. of the standard representation of
  ${\mathfrak {su}}(n,1)$ or of its contragredient representation).
\end{thm}

In this theorem, Raghunathan used Matsushima-Murakami's result where
$L^2$-cohomology is used. We observe that as long as we use
$L^2$-cohomology, this theorem still holds for non-uniform lattices.
This issue will be dealt with in section \ref{nonuniform}.

\subsection{Standard representation of $\mathfrak{sp}(4,\bc)$}
\label{symplectic}

In the previous section, we used the symplectic form with respect to
the standard basis of $\bc^4$
$$ Q=\left[  \begin{matrix}
    0   &  \begin{matrix}
           1 & 0 \\
           0 & -1 \end{matrix}   \\
     \begin{matrix}
           -1 & 0 \\
           0 & 1 \end{matrix}    & 0  \end{matrix}
           \right].$$
Then the Lie algebra $\mathfrak{sp}(4,\bc)$ consists of complex
matrices $\displaystyle \left[ \begin{matrix}
        A & B \\
        C & D \end{matrix} \right]$ such that
$$A^t \left[\begin{matrix}
           1 & 0 \\
           0 & -1 \end{matrix}\right] + \left[\begin{matrix}
           1 & 0 \\
           0 & -1 \end{matrix}\right] D=0,$$
$$    C^t \left[\begin{matrix}
           -1 & 0 \\
           0 & 1 \end{matrix}\right] + \left[\begin{matrix}
           1 & 0 \\
           0 & -1 \end{matrix}\right] C=0,$$
$$  B^t \left[\begin{matrix}
           1 & 0 \\
           0 & -1 \end{matrix}\right] + \left[\begin{matrix}
           -1 & 0 \\
           0 & 1 \end{matrix}\right] B=0.$$
Then an obvious choice of a Cartan subalgebra $\mathfrak{h}$ is
$$\left[ \begin{matrix}
       \begin{matrix}
       x & 0 \\
       0 & y \end{matrix}  &  0 \\
       0       &   \begin{matrix}
                    -x & 0 \\
                    0  & -y \end{matrix} \end{matrix}\right].$$

Let $L_1$ and $L_2 \in\mathfrak{h}^*$ be defined by $L_1 (x,y)=x$,
$L_2(x,y)=y$. Then the natural action of $\mathfrak{sp}(4,\bc)$ on
$\bc^4$ has the four standard basis vectors $e_1,e_2,e_3,e_4$ as
eigenvectors with weights $L_1,L_2,-L_1,-L_2$. The highest weight is
$L_1$.

\subsection{Representation of $\mathfrak{so}(5,\bc)$.}

We shall use the isomorphism of $\mathfrak{sp}(4,\bc)$ to
$\mathfrak{so}(5,\bc)$. It arises from the following geometric
construction.

Let $V=\bc^4$ and $\omega$ be the symplectic form defined as before.
Then $$\wedge^2 V^* \otimes \wedge^2 V^* \ra \bc$$
$$ \alpha \otimes \beta \ra \frac{\alpha \wedge \beta}{\omega \wedge \omega},$$
is a nondegenerate quadratic form $P$ on $\wedge^2 V^*$. Here since
both $\alpha\wedge \beta$ and $\omega\wedge \omega$ are 4-forms,
there is a constant $c$ so that $\alpha\wedge \beta=c\omega\wedge
\omega$, so the quotient should be understood as such a constant.
Take the orthogonal complement $W$ of $\bc \omega$ with respect to
this quadratic form. Any matrix $A$ acts on 2-forms as follows:
$A\alpha(v,w)=\alpha(Av,Aw)$. Then $Sp(4,\bc)$ leaves $W$ invariant
and acts orthogonally on it. This gives a map from $Sp(4,\bc)$ to
$SO(5,\bc)=SO(W)$, which turns out to be an isomorphism.

Next, we relate the choice of Cartan subalgebra for
$\mathfrak{sp}(4,\bc)$ made in the preceding paragraph to the
standard choice for $\mathfrak{so}(5,\bc)$.

We first compute the Lie algebra isomorphism derived from the group
isomorphism.

Let $z_1,z_2,z_3,z_4$ be standard coordinates of $\bc^4$ so that
$dz_1\wedge dz_3+ dz_4 \wedge dz_2=\omega$. Let $\omega_6=\omega$
and $$\omega_5=dz_1\wedge dz_2+ dz_3\wedge dz_4,$$
$$\omega_4=dz_1\wedge dz_4 + dz_2 \wedge dz_3,$$
$$\omega_1=i(dz_1\wedge dz_4 - dz_2\wedge dz_3),$$
$$\omega_2=i(dz_1\wedge dz_2- dz_3\wedge dz_4),$$
$$\omega_3=i(dz_1\wedge dz_3-  dz_4 \wedge dz_2).$$
This is an orthonormal basis of $\wedge^2 V^*$.

Let $A_t\in Sp(4,\bc)$ so that $A_0=I$ and $\frac{d}{dt}|_{t=0}
A_t=X\in \mathfrak{sp}(4,\bc)$. Then for one-forms $\alpha,\beta$,
one can figure out the action of $X$ on two-forms to see that $X
(\alpha\otimes \beta)= \frac{d}{dt}|_{t=0} A_t (\alpha\otimes
\beta)=(X\alpha)\otimes \beta + \alpha\otimes (X\beta)$. Then
$$X(\alpha \wedge \beta)=(X\alpha)\wedge \beta + \alpha \wedge
(X\beta).$$

To make computation easier, we choose a basis of $W$ as
$$v_1=\frac{\omega_1+ i \omega_4}{\sqrt 2},$$
$$v_3=\frac{\omega_1- i \omega_4}{\sqrt 2},$$
$$v_2=\frac{\omega_2+i\omega_5}{\sqrt 2},$$
$$v_4=\frac{\omega_2-i\omega_5}{\sqrt 2},\ v_5=\omega_3.$$ With
respect to this basis, the symmetric bilinear form $P$ has
$P(v_1,v_3)=1=P(v_2,v_4)=P(v_5,v_5)$ and $P(v_i,v_j)=0$ for all
other pairs. With respect to this $P$, one can easily see that a
Cartan subalgebra of $\mathfrak{so}(5,\bc)=\mathfrak{so}(W;P)$ can
be chosen as the set of matrices of the form
$$\left [ \begin{matrix}
  x & 0 & 0 & 0 & 0 \\
  0 & y & 0 & 0 & 0 \\
  0 & 0 & -x& 0 & 0 \\
  0 & 0 & 0 & -y& 0 \\
  0 & 0 & 0 & 0 & 0 \end{matrix}\right].$$

 Let $(x,y,z,w)$ denote a diagonal matrix in $\mathfrak{sp}(4,\bc)$. Then
one can easily compute that $$(1,0,-1,0) v_1=v_1,
(1,0,-1,0)v_3=-v_3,$$
$$(1,0,-1,0)v_2=v_2,(1,0,-1,0)v_4=-v_4,(1,0,-1,0)v_5=0.$$
Similarly
$$(0,1,0,-1)v_1=-v_1,(0,1,0,-1)v_3=v_3,$$
$$(0,1,0,-1)v_2=v_2, (0,1,0,-1)v_4=-v_4,(0,1,0,-1)v_5=0.$$

 So the element, $\left[ \begin{matrix}
  1 &0&0&0 \\
  0 &0&0&0 \\
  0 &0&-1&0 \\
  0&0&0&0 \end{matrix}\right]$, in a Cartan subalgebra of
  $\mathfrak{sp}(4,\bc)$ corresponds to an element in a Cartan
  subalgebra of $\mathfrak{so}(5,\bc)$,
  $$\mathfrak{h_1}=\left[ \begin{matrix}
  1 & 0 & 0 & 0 & 0 \\
  0 & 1 & 0 & 0 & 0 \\
  0 & 0 & -1& 0 & 0 \\
  0 & 0 & 0 & -1 &0 \\
  0 & 0 & 0 & 0 & 0 \end{matrix} \right].$$
Similarly  $\left[ \begin{matrix}
  0 &0&0&0 \\
  0 &1&0&0 \\
  0 &0&0 &0 \\
  0&0&0&-1 \end{matrix}\right]$ corresponds to
$$\mathfrak{h_2}=\left[ \begin{matrix}
  -1 & 0 & 0 & 0 & 0 \\
  0 & 1 & 0 & 0 & 0 \\
  0 & 0 & 1& 0 & 0 \\
  0 & 0 & 0 & -1 &0 \\
  0 & 0 & 0 & 0 & 0 \end{matrix} \right].$$



This representation under the isomorphism to $\mathfrak{sp}(4,\bc)$
is different from the standard representation of
$\mathfrak{so}(5,\bc)$ on $\bc^5$ as we will see below.

\begin{lemma}
\label{spin} The highest weight of the standard
  representation
of $\mathfrak{so}(5,\bc)$ on $\bc^5$ is not a multiple of the
highest weight of
 the representation coming from
$\mathfrak{sp}(4,\bc)$ on $\bc^4$.
\end{lemma}

\begin{pf}
With respect to the symmetric bilinear form $P$ as before, a Cartan
subalgebra of $\mathfrak{so}(5,\bc)$
 is the set of diagonal matrices
$(x,y,-x,-y,0)$ as noted above. Then the standard representation of
$\mathfrak{so}(5,\bc)$ on $\bc^5$ has eigenvectors, the standard
basis $e_1,e_2,e_3,e_4,e_5$, with eigenvalues $L_1,L_2,-L_1,-L_2,0$.
This has the highest weight $L_1$.

The standard representation of $\mathfrak{sp}(4,\bc)$ on $\bc^4$ has
the highest weight $L_1$ as we saw in the previous section. Note
that the Cartan subalgebra of $\mathfrak{sp}(4,\bc)$ is generated by
the diagonal matrices $(1,0,-1,0)$ and $(0,1,0,-1)$ with dual basis
$L_1$ and $L_2$. Then under the isomorphism from
$\mathfrak{sp}(4,\bc)$ to $\mathfrak{so}(\bc \omega^\perp)$, these
two diagonal matrices are mapped to diagonal matrices
$\mathfrak{h}_1=(1,1,-1,-1,0)$ and $\mathfrak{h}_2=(-1,1,1,-1,0)$.
Let $L_1',L_2'$ be the images of $L_1,L_2$ under this isomorphism.
Then in terms of the standard dual basis $L_1,L_2$ of the Cartan
subalgebra of $\mathfrak{so}(5,\bc)$,
$$L_1'=\frac{L_1 + L_2}{2}, L_2'=\frac{L_2-L_1}{2}.$$
So the representation coming from the standard representation of
$\mathfrak{sp}(4,\bc)$ on $\bc^4$ has highest weight
$\frac{L_1+L_2}{2}$. Actually this is the highest weight of the spin
representation.
\end{pf}

\begin{co}
\label{rag} Let $\Gamma \subset Sp(1,1)$ be a uniform lattice. Then
$H^1(\Gamma, \bh^2)=0$ where $\bh^2$ is denotes the standard
representation of $Sp(1,1)$ restricted to $\Gamma$.
\end{co}

\begin{pf}
View $\bh^2$ as $\bc^4$ with $Sp(1,1)$ acting on it. If we
complexify the real Lie algebra $\mathfrak{sp}(1,1)$, we get
$\mathfrak{sp}(4,\bc)$. Since the standard representation of
$\mathfrak{sp}(4,\bc)$ on $\bc^4$ is different from the standard
representation of $\mathfrak{so}(5,\bc)$ on $\bc^5$ with highest
weight $L_1$, Theorem \ref{raghunathan} (Theorem 1 of Raghunathan
\cite{Ra}) applies, and $H^1(\Gamma, \bh^2)=0$.
\end{pf}

\section{Proof of Theorem \ref{4} (uniform case)}
\label{4u}

Let $\Gamma \subset Sp(1,1)$ be a uniform lattice. Denote by $\rho$
the embedding $\Gamma \ra Sp(1,1) \ra Sp(2,1)$. Let $G=Sp(2,1)$,
$H=\left [\begin{matrix}
Sp(1) & 0\\
0  & Sp(1,1)\end{matrix}\right]\subset G$. As was seen in section
\ref{pre}, the adjoint representation of $G$ restricted to $H$
splits as a direct sum $\mathfrak{sp}(2,1)=\mathfrak{sp}(1)\oplus
\bh^2 \oplus\mathfrak{sp}(1,1)$, thus
$\mathfrak{g}/\mathfrak{h}=\bh^2$, restricted to $Sp(1,1)$, is the
standard representation of $Sp(1,1)$. Corollary \ref{rag} asserts
that $H^1(\Gamma,\bh^2)$ vanishes. Therefore
$H^1(\Gamma,\mathfrak{g}_{\rho}/\mathfrak{h}_{\rho})=0$. According
to Proposition \ref{relative}, this implies that homomorphisms
$\Gamma\ra Sp(2,1)$ which are close enough to $\rho$ can be
conjugated into $H$, i.e. leave a quaternionic line invariant.

Since the subgroup of the form
$$\left[ \begin{matrix}
    Sp(1) & 0 & 0 \\
    0     & I & 0 \\
    0    & 0 & Sp(1,1) \end{matrix} \right]$$
stabilizes a quaternionic line $(0,0,\cdots,0,\bh)$ in the ball
model, we obtain
\begin{co}Let $\Gamma \subset Sp(1,1)$ be a uniform lattice.
Embed $\Gamma$ into $Sp(n,1)$ as a subgroup which stabilizes a
quaternionic line. Then every small deformation of $\Gamma$ in
$Sp(n,1)$ stabilizes a quaternionic line.
\end{co}

\section{3-manifold case}
\label{low}

In this section, we prove Theorem \ref{3} for uniform 3-dimensional
hyperbolic lattices. Let $\Gamma\subset Spin(3,1)^0$ be a uniform
lattice. According to Proposition \ref{relative}, local deformations
of the standard representation $\rho_0:\Gamma \ra Spin(3,1)^0 \ra
Spin(4,1)^0 =Sp(1,1)\ra Sp(2,1)$ which do not stabilize a
quaternionic line, are encoded in $H^1(\Gamma,\bh^2)$. We want to
show that this first cohomology is zero. The complexified Lie
algebra of $SO(3,1)$ is $\mathfrak{so}(4,\bc)$. In the notations of
section \ref{rep}, the symmetric bilinear form $P$ has a basis
$v_1,v_2,v_3,v_4$ so that $P(v_1,v_3)=P(v_2,v_4)=1$ and
$P(v_i,v_j)=0$ for all other pairs. The Cartan subalgebra of
$\mathfrak{so}(4,\bc)$ is the set of diagonal matrices
$(x,y,-x,-y)$. Then as in Lemma \ref{spin}, the standard
representation of $\mathfrak{so}(4,\bc)$ on $\bc^4$ has a character
which is not a multiple of the character of the representation
coming from $\mathfrak{so}(4,\bc)\subset \mathfrak{sp}(4,\bc)$. Then
by Raghunathan's theorem \ref{raghunathan}, $H^1(\Gamma,\bh^2)=0$.
Proposition \ref{relative} ensures that neighboring homomorphisms
$\Gamma \ra Sp(2,1)$ stabilize a quaternionic line.

\section{Non-uniform lattices}
\label{nonuniform}

We used Raghunathan's theorem \cite{Ra} to prove our main theorem
when $\Gamma$ is a uniform lattice. In this section we discuss how
it generalizes, with restrictions, to nonuniform lattices.

The key point is whether Matsushima-Murakami's vanishing theorem
that Raghunathan used still holds in non-uniform case. To apply
Matsushima-Murakami's theorem, one has to use $L^2$-cohomology.

Recall that under the subgroup $\left [\begin{matrix}
Sp(1) & 0\\
0  & Sp(1,1)\end{matrix}\right]$, the adjoint representation of
$Sp(2,1)$ splits as a direct sum
$\mathfrak{sp}(2,1)=\mathfrak{sp}(1)\oplus \bh^2
\oplus\mathfrak{sp}(1,1)$. Let $\rho$ denote the representation of
$Sp(1,1)$ corresponding to the $\bh^2$ summand. Let
$M=H^4_\br/\Gamma$ be a finite volume manifold. View $\Gamma$ as a
subgroup of $Sp(1,1)$, denote by $\rho_0$ the restriction of $\rho$
to $\Gamma$. Let $E$ be the associated flat bundle over $M$ with
fibre $\bh^2$. It is well-known that
$$H^1(\Gamma,\rho_0)=H^1_{dR}(M,E)$$
where $H^1_{dR}(M,E)$ is de Rham cohomology of smooth $E$-valued
differential forms over $M$. We will denote this de Rham cohomology
by $H^1 (M,E)$.

In Matsushima-Murakami's proof, specific metrics on fibres of $E$,
depending on base points, are used. More precisely, fix a maximal
compact subgroup $K$ of $Sp(1,1)$. Let
$\mathfrak{sp}(1,1)=\mathfrak{t}\oplus\mathfrak{p}$ be the
corresponding Cartan decomposition. Fix a positive definite metric
$\langle\ ,\ \rangle_F$ on $\bh^2$ so that $\rho(K)$ is unitary and
$\rho(\mathfrak{p})$ is hermitian symmetric. Then, for two elements
$v$, $w$ in the fibre over a point $g\in G$, one defines
$$\langle v, w \rangle = \langle \rho(g)^{-1}v, \rho(g)^{-1}w \rangle_F.$$
Here is a concrete construction of such a metric on $\bh^2$. As
before, $\bh^{1,1}=\bh^2$ is equipped with the signature
$(1,1)$-metric
$$Q=|q_1|^2-|q_2|^2.$$
Then for each negative $\bh$-line $L$ in $\bh^{1,1}$, there exists a
positive definite $\bh$-Hermitian metric defined by $-Q|_L \oplus
Q|_{L^\perp}$ where $L^\perp$ is the orthogonal complement of $L$
with respect to $Q$.

A unit speed ray in $H^4_\br=H^1_\bh$ in terms of $\bh^{1,1}$
coordinates, can be written as $l_t=\{q_1=\delta_t q_2\}$ where
$\delta_t=\frac{e^t-1}{e^t+1},\ 0\leq t \leq \infty$. Note that here
we normalize the metric so that its sectional curvature is $-1$.
This can be easily computed considering a unit speed ray $r(t)$ in a
ball model emanating from the origin, and $r(t)$ corresponds to the
point $(r(t),1)$ in the hyperboloid model.

Now we want to know how the metric varies along $l_t$ as $t \ra
\infty$. Let $v=(v_1,v_2)\in \bh^{1,1}$. It is easy to see that
$$b_t=(\frac{1}{\sqrt{1-\delta_t^2}},\frac{\delta_t}{\sqrt{1-\delta_t^2}}),$$
$$a_t=(\frac{\delta_t}{\sqrt{1-\delta_t^2}},\frac{1}{\sqrt{1-\delta_t^2}})$$
are unit vectors on $l_t^\perp,l_t$ respectively. Then $l_t$
component of $v$ is
$$(\frac{\delta_t v_2 +\delta^2_t
  v_1}{1-\delta_t^2},\frac{v_2+\delta_t v_1}{1-\delta_t^2})$$
and $l_t^\perp$ component is
$$(\frac{\delta_t v_2 +  v_1}{1-\delta_t^2},\frac{\delta_t^2 v_2 +
\delta_t v_1}{1-\delta_t^2}).$$ Then it is easy to calculate the
square of the length of $v$ on $l_t$, which is
$$\frac{1+\delta_t^2}{1-\delta_t^2}
[|v_1|^2 + |v_2|^2]+ 2 \frac{\delta_t}{1-\delta_t^2}
(v_1\overline{v_2}+v_2\overline{v_1})$$
$$=\frac{2\delta_t}{1-\delta_t^2}|v_1+v_2|^2+\frac{1-\delta_t}{1+\delta_t}
(|v_1|^2+|v_2|^2).$$ In conclusion, the square of the length of $v$
grows like $e^t |v|^2$ along the ray $l_t$ in general. But for
$v_1+v_2=0$, it grows like $ e^{-t}|v|^2$ along the ray. This is the
case when the deformation consists in parabolic elements fixing a
point $(0,1)$ (in the ball model) and not stabilizing the
quaternionic line $\{(0,\bh)\}$. See Lemma \ref{stab}. These
estimates will be used  below.

Let $M=M_{\geq \epsilon}\cup M_{\leq \eta}$ be the thick-thin
decomposition of $M$ so that $\eta> \epsilon$ and $M_{\leq \eta}$ is
a standard cusp part of $M$. Assume for simplicity that the cuspidal
part is connected. It is well-known that $M_{\leq \eta}$ is
homeomorphic to $T\times \br^+$ with $ds^2=e^{-2r}ds^2_T + dr^2$
where $T$ is a flat closed 3-manifold, $r$ denotes distance from
$T\times\{0\}$, and $M_{\geq \epsilon}\cap M_{\leq \eta}$ is
$T\times [0,1]$.

%
%

Let $\pi:T\times \br^+ \ra T$ be the projection on the first factor.
Since $H^k(T)=H^k(M_{\leq \eta})$ by $\pi^*$, we want to show that
$L^2H^k(M_{\leq \eta})=H^k(T)$, to show that $H^k(M_{\leq
\eta})=L^2H^k(M_{\leq \eta})$. Let $\alpha$ be a $k$-form on $T$.
Then $|\pi^*\alpha|\sim e^{\frac{r}{2}}|\alpha|e^{kr}$ where $r$ is
the distance from the boundary of the thin part. Here
$e^{\frac{r}{2}}$ comes from the fibre metric and $e^{kr}$ comes
from the base metric. Then
$$||\pi^*\alpha ||^2_{L^2}=\int |\alpha|^2
e^{2kr+r}e^{-3r}ds_T dr\leq ||\alpha||^2_{L^2(T)}\times C <\infty$$
if $2k+1<3$. So the pull-back form $\pi^*\alpha$ is always a
$L^2$-form on $M_{\leq \eta}$ if $\alpha$ is a 0-form.

So we obtained

\begin{lemma}
\label{l2} For a finite volume real 4-dimensional hyperbolic
manifold $M$, $H^0 (M_{\leq \eta},E)=L^2 H^0 (M_{\leq \eta},E)$.
\end{lemma}

\begin{pf}
For any $\alpha \in H^*(T,E)=H^*(M_{\leq \eta},E)$, its pull-back
$\pi^*\alpha$ is a $L^2$-form on $M_{\leq \eta}$ for $*=0$ as noted
above. So any element in $H^0(M_{\leq \eta},E)$ has an
$L^2$-representative.
\end{pf}

Unfortunately, we cannot conclude that $H^1(M,E)=L^2 H^1(M,E)$. This
hinders us from generalizing our theorem to non-uniform lattices.
Our generalization involves a restriction on the representation.

\begin{Prop}
Let $M$ be a finite volume hyperbolic 3-manifold so that
$M=H^3_\br/\Gamma$. Then all small deformations of $\Gamma\subset
SO(3,1)\subset Sp(1,1)$ preserving parabolicity still stabilizes a
quaternionic line. The same thing holds for a finite volume
hyperbolic 4-manifold.
\end{Prop}

\begin{pf}
We give a proof only in dimension 3, since the 4-dimensional case
can be obtained by the same method. Since $M$ has finite volume, its
boundary consists of tori $T_i$. Let $\rho_0 :\pi_1(M)\ra
Spin(3,1)^0\subset Sp(1,1)\subset Sp(2,1)$ be a natural
representation.

If $\rho_t (\pi_1(\partial M))$ is parabolic for all small $t$, by
Lemma \ref{stab}, it can contribute to the $\bh^2$ summand of
${\mathfrak{sp}(2,1)}$.
But in this case, it can be represented by an $L^2$ form. The
argument goes briefly as follows.

Let $\rho_t:\pi_1(M)\rightarrow Sp(2,1)$ be an one-parameter family
of deformations so that $\rho_t(\pi_1(\partial M))$ is all
parabolic. Let $N$ be the $\epsilon$-thick part of $M$. Then
$\partial N$ consists of tori and the universal cover of it in
$H^3_\br$ are horospheres.
%
%
Fix a component of $\tilde {\partial N}$ which is a horosphere $H$
corresponding to a component $T$ of $\partial N$.
 Conjugating $\rho_t$ by $g_t$ which depend smoothly on $t$ if necessary, we may assume
that $\rho_t(\pi_1(T))$ leaves invariant a common horosphere $H'$ in
$H^2_\bh$. Such a choice of $g_t$ is possible by the following
argument. Let $a$ be an element in $\pi_1(T)$ such that all
$\rho_t(a)$ are parabolic. The subset $P$ of $Sp(2,1)$ consisting of
parabolic elements is a smooth manifold at $\rho_0(a)$, and the map
from $P$ to $\partial{H^2_\bh}$ associating to each element in $P$
its unique fixed point is smooth in a neighborhood of $\rho_0(a)$.

We may assume that $H'$ is based at $(0,1)$ (in the ball model).
Then by Lemma \ref{stab}, the contribution of this deformation to
the $\bh^2$ summand is contained in the subset
$\{(x,y)|x+y=0\}\subset \bh^2$. This will help us out.

Let $\omega$ be a differential form representing the infinitesimal
deformation $\frac{d}{dt}\rho_t$ on this cusp. Since
$\rho_t(\pi_1(T))$ fixes $(0,1)$, $\omega$ takes its values in the
subalgebra $\mathfrak{s}\subset \mathfrak{sp}(2,1)$ of Killing
fields on $H^2_\bh$ which vanish at $(0,1)$ and which are tangent to
the horospheres centered at $(0,1)$. Therefore the norm of vectors
of $\mathfrak{s}$ decays along a geodesic pointing to $(0,1)$, at
speed controlled by the maximal sectional curvature (in our case,
which is the direction away from a quaternionic line,
$-\frac{1}{4}$). In our situation, we are only concerned with the
subspace $\{(v_1,v_2)|v_1+v_2=0\}\subset \bh^2$. So along the ray
the squared norm decays like $e^{-r}|v|^2$ asymptotically.

Then integrating along a geodesic ray, we see that the 1-form
$\omega$
 defined on the cusp is in $L^2$ on the cusp.
In more details, let the cusp be $T\times [0,\infty)$ with
coordinates $(x,y,r)$, and the metric $ds^2=e^{-2r}ds_T^2 + dr^2$,
then the volume form on this cusp is  $e^{-2r}dS_Tdr$. Note that we
take a metric on $H^3_\br$ whose sectional curvature is $-1$. Then
along $[0,\infty)$, the orthonormal basis is $\{e^r
\frac{\partial}{\partial x}, e^r \frac{\partial}{\partial y},
\frac{\partial}{\partial r}\}$. Then at $(x,y,r)$, the norm of
$\omega$ is
$$|\omega(e^r \frac{\partial}{\partial x})|^2+
|\omega(e^r \frac{\partial}{\partial y})|^2$$ since
$\omega(\frac{\partial}{\partial r})=0$.

So
$$\int_{T\times [0,\infty)} ||\omega||^2 dVol=
\int_0^\infty e^{-r}e^{2r} e^{-2r}\int _T ||\omega_T||^2 dS_T
dr<\infty$$ where $e^{-r}$ comes from the norm decay on
$\{(v_1,v_2)|v_1+v_2=0\}$, $e^{2r}$ comes from the decay of the
metric on $H^3_\br$ along the ray (one should take an orthonormal
basis $\{e^r \frac{\partial}{\partial x}, e^r
\frac{\partial}{\partial y}, \frac{\partial}{\partial r}\}$ along
the ray).

 We do this for each cusp of $M$. Let $\omega_i$ be a 1-form which is a $L^2$-representative of the deformation
 $\frac{d}{dt}\rho_t$ on the $i$-th cusp of $M$. Let $\alpha$ be a global 1-form representing the deformation $\frac{d}{dt}\rho_t$. Then
$$\omega_i=\alpha + d\phi_i$$
where $\phi_i$ is  a function defined on the $i$-th cusp. Let $\phi$
be the union of $\phi_i$ and $\xi$ be a smooth function so that
$\xi=1$ on cusps and $0$ outside cusps. Let
\begin{eqnarray*}
\omega'&=&\alpha + d(\xi \phi)\\
&=&\alpha +\phi d\xi + \xi d\phi.
\end{eqnarray*}
Then on each cusp, $\omega'=\alpha + d\phi_i=\omega_i$. Thus
$\omega'$ is in $L^2$ and $[\omega']=[\alpha]$.

Now again we can use Matsushima-Murakami's result for this case. See
\cite{KM1, KM2} for a similar argument in complex hyperbolic space.

So we proved the theorem.
%
\end{pf}

We wonder whether the theorem holds without the assumption of
preserving parabolicity.



\section{Bending representations}
\label{bend}

Let $G$ be an algebraic group. The Zariski closure of a subgroup $H$
of $G(\br)$ is denoted by $\bar{H}$.

Let $X$ be a compact orientable hyperbolic $n$-manifold which splits
into two submanifolds with totally geodesic boundary $V$ and $W$,
exchanged by an involution that fixes their common boundary. Such
manifolds exist in all dimensions, \cite{Mill}. Then $\Gamma=\pi_1
(X)$ splits as an amalgamated sum $\Gamma=A\star_{C}B$ where
$A=\pi_1 (V)$, $B=\pi_1 (W)$ and $C=\pi_1 (\partial V)$. Here,
$\bar{A}=\bar{B}=PO(n,1)^0$ and $\bar{C}=PO(n-1,1)^0$.

Now embed $PO(n,1)^0$ into a larger group $G$. Let $c$ belong to the
centralizer $Z_{G}(C)$. Consider the subgroup $\Gamma_c =A\star_C
cBc^{-1}$. When $c$ is chosen along a curve in $Z_{G}(C)$, one
obtains a special case of W. Thurston's {\em bending deformation},
\cite{Th} chapter 6. In this section, we analyze the Zariski closure
of $\Gamma_c$ in case $G=PSp(m,1)$ is the isometry group of
$m$-dimensional quaternionic hyperbolic space, $m\geq n$ and
$PO(n,1)^0 \ra PSp(n,1)\ra PSp(m,1)$ in the obvious manner.

\subsection{The first bending step}

We find it convenient to use a geometric language, and establish a
dictionary between subgroups of $G=PSp(m,1)$ and totally geodesic
subspaces of $X=H_{\bh}^{m}$.

\begin{lemma}
\label{normalizer} The subgroup of $G$ that leaves
$Y=H_{\br}^{n}\subset X$ invariant is the normalizer of
$H=PO(n,1)^0$ in $G$.
\end{lemma}

\begin{pf}
If $aHa^{-1}=H$, $a$ maps the orbit $Y$ of $H$ to itself.
Conversely,  $Y$ is the only orbit of $H$ in $X$ which is totally
geodesic. If $a\in G$ normalizes $H$, then $a$ maps $Y$ to itself.
\end{pf}

Second, let us determine the space of available parameters for
bending, i.e. elements which commute with $C$.

\begin{lemma}
\label{centralizer} Let $m\geq n\geq 2$. Let $L=PO(n-1,1)^0 \subset
PO(n,1)^0 \subset PSp(n,1)\subset PSp(m,1)=G$. Let $C\subset L$ be a
Zariski dense subgroup. Then the centralizer $Z_{G}(C)$ consists of
isometries which fix $P=H_{\br}^{n-1}$ pointwise. As a matrix group,
$Z_{G}(C)=Sp(m-n+1)Sp(1)$.
\end{lemma}

\begin{pf}
Clearly, $Z_{G}(C)=Z_{G}(L)$. $L$ stabilizes the totally geodesic
subspace $P=H_{\br}^{n-1}$ of the symmetric space $X=H_{\bh}^{m}$ of
$G$. If $a\in G$ centralizes $L$, then $a$ normalizes it, thus it
maps $P$ to itself, by Lemma \ref{normalizer}. Furthermore, the
restriction of $a$ to $P$ belongs to the center of $Isom(P)= L$,
thus is trivial. In other words, $a$ fixes each point of $P$.
Conversely, isometries of $X$ which fix every point of $P$
centralize $L$ and thus $C$. Indeed, $L$ is generated by geodesic
symmetries with respect to points of $P$, and these commute with
isometries fixing $P$. To get the matrix expression of $Z_{G}(C)$,
view $X$ as a subset of quaternionic projective $m$-space. Then for
every vector $y\in\br^n$, extended with zero entries to give a
vector in $\br^{m+1}$, there exists a quaternion $q(y)$ such that
$a(y)=yq(y)$. This implies that $a$ lifted as a matrix in $Sp(m,1)$
is block diagonal,
\begin{eqnarray*}
a=\left[
\begin{matrix}
     qI_n & 0\\
     0 & D \end{matrix} \right],
\end{eqnarray*}
with blocks of sizes $n$ and $m-n+1$ respectively, $q\in Sp(1)$ and
$D\in Sp(m-n+1)$. This product group maps to a subgroup of
$PSp(m,1)$ which is traditionnally denoted by $Sp(m-n+1)Sp(1)$.
\end{pf}

The dictionary continues with a correspondance between Zariski
closures in simple groups and totally geodesic hulls in symmetric
spaces.

\begin{lemma}
\label{zariski} Let $Y_1 ,\ldots,Y_k$ be totally geodesic subspaces
of a symmetric space $X$. Then $Isom(Y_j)$ naturally embeds into
$G=Isom(X)$. Furthermore, the Zariski closure of $\bigcup_j
Isom(Y_j)$ equals $Isom(Z)$ where $Z$ is the smallest totally
geodesic subspace of $X$ containing $\bigcup_j Y_j$.
\end{lemma}

\begin{pf}
For $x\in X$, let $\iota_x$ denote the geodesic symmetry through
$x$. Since $X$ is symmetric, $\iota_x$ is an isometry. Such
involutions generate $Isom(X)$. If $Y\subset X$ is totally geodesic,
then $Y$ is invariant under all $\iota_y$, $y\in Y$. Therefore $Y$
is again a symmetric space, with isometry group generated by the
restrictions to $Y$ of the $\iota_y$. In particular, the subgroup of
$G$ generated by the $\iota_y$, $y\in Y$, is isomorphic to
$Isom(Y)$.

If $\gamma$ is a geodesic joining points $x\in Y_i$ and $y\in Y_j$,
then $\iota_x$ and $\iota_y$ leave $\gamma$ invariant. Their
restrictions to $\gamma$ generate an infinite dyadic group. The
Zariski closure of this group contains all $\iota_z$ where
$z\in\gamma$. Therefore the Zariski closure of $Isom(Y_i)\cup
Isom(Y_j)$ contains $\iota_z$ for all $z$ belonging to the union of
all geodesics intersecting both $Y_i$ and $Y_j$. Since the totally
geodesic closure $Z$ is obtained by iterating this operation, one
concludes that the Zariski closure of $\bigcup_j Isom(Y_j)$ contains
$Isom(Z)$. Conversely, since $Isom(Z)$ is an algebraic subgroup in
$G$, it is contained in the Zariski closure.
\end{pf}

\begin{lemma}
\label{totally} Let $Y=H_{\br}^{n}\subset H_{\bh}^{n}=X$. Let $Z$ be
a totally geodesic subspace of $X$ such that $Y\subsetneq Z
\subsetneq X$. Assume that $Z$ contains $a(Y)$ where $a\in G$ fixes
pointwise a hyperplane $P$ of $Y$ but does not leave $Y$ invariant.
Then there is an isometry of $X$ fixing $Y$ pointwise and mapping
$Z$ to $H_{\bc}^{n}$.
\end{lemma}

\begin{pf}
View the restriction of $TX$ to $Y$ as a vector bundle with
connection $\nabla$ on $Y$. Then $TZ_{|Y}$ is a parallel subbundle,
therefore, for $y\in Y$, $T_y Z$ is invariant under the holonomy
representation $Hol(\nabla,y)$, which we now describe.

View $Y$ as a sheet of the hyperboloid in $\br^{n+1}$. Then a point
$y$ represents a unit vector, still denoted by $y$, in $\br^{n+1}$.
View $X$ as a subset of quaternionic projective space. Then the
point $y$ also represents the quaternionic line $\bh y$ it
generates. Such lines form the tautological quaternionic line bundle
$\tau$ over $X$, a subbundle of the trivial bundle $\bh^{n+1}$
equipped with the orthogonally projected connection. As a connected
vector bundle, $TX=Hom_{\bh}(\tau,\tau^{\bot})$. When restricted to
$Y$, $\tau$ comes with the parallel section $y$. Therefore
$TX_{|Y}=\tau^{\bot}=TY\otimes\bh$. In other words, $TX_{|Y}$ splits
as a direct sum of 4 parallel subbundles, each of which is
isomorphic to $TY$. It follows that $Hol(\nabla,y)$ is the direct
sum of four copies of the holonomy of the tangent connection, which
is the full special orthogonal group $SO(n)$. One of these copies is
$T_y Y$, the other are its images under an orthonormal basis
$(I,J,K)$ of imaginary quaternions acting on the right.

Let us show that $Z$ contains a copy of $H_{\bc}^{n}$. Let $a\in G$
fix a hyperplane $P\subset Y$ pointwise. According to Lemma
\ref{centralizer}, $Fix(P)=Sp(1)Sp(1)$, so $a$ is given by two unit
quaternions $q$ and $d$. Pick an origin $y\in P$. Let $u\in T_y Y$
be a unit vector orthogonal to $P$. On $T_y X=T_y Y\otimes\bh$, $a$
acts by the identity on $T_y P$ and maps $u$ to $duq^{-1}$. Since
$u$ is a real vector, $a(u)=udq^{-1}\in T_y Y\oplus (T_y Y)i$ where
$i=\Im m(dq^{-1})$. Up to conjugating by an element of the $Sp(1)$
subgroup of $G$ that fixes $Y$ pointwise, one can assume that $i$ is
proportional to $I$, i.e. $T_y Z$ contains $uI$. By assumption,
$uI\notin T_y Y$. By $SO(n)$ invariance, $T_y Z$ contains $T_y
Y\oplus (T_y Y)I=T_y H_{\bc}^{n}$, therefore $Z$ contains
$Y'=H_{\bc}^{n}$.

Now $TZ_{|Y'}$ is a parallel subbundle of $TX_{|Y'}$, thus $T_y Z$
is $U(n)$-invariant. Under $U(n)$, $T_y X$ splits into only 2
summands. Since $Z\not=X$, $T_y Z=T_y Y'$, i.e. $Z=Y'$.
\end{pf}

Along the way, we proved the following.

\begin{lemma}
\label{totallycomplex} Let $Y'=H_{\bc}^{n}\subset H_{\bh}^{n}=X$.
Let $Z$ be a totally geodesic subspace of $X$ containing $Y'$. Then
either $Z=X$ or $Z=Y'$.
\end{lemma}

\begin{co}
\label{complex} After bending in $PSp(n,1)$, a Zariski dense
subgroup of $PO(n,1)^0$ becomes Zariski dense in a conjugate of
$PU(n,1)$.
\end{co}

\begin{pf}
Let $\Gamma=A\star_C B$ be Zariski dense in $PO(n,1)^0$, with $C$
Zariski dense in $PO(n-1,1)^0$. In other words, $\Gamma$ leaves
$Y=H_{\br}^{n}$ invariant, and $C$ leaves $P=H_{\br}^{n-1}$
invariant. Lemma \ref{centralizer} allows to select an $a\in Z_G
(C)$ which does not map $Y$ to itself. Lemma \ref{totally} shows
that the smallest totally geodesic subspace of $X=H_{\bh}^{n}$
containing $Y$ and $a(Y)$ is congruent to $H_{\bc}^{n}$. According
to Lemma \ref{zariski}, this means that the bent subgroup $A\star_C
aBa^{-1}$ is Zariski dense in a conjugate of $PU(n,1)$.
\end{pf}

Therefore, to obtain a Zariski dense subgroup in $PSp(m,1)$, $m\geq
n$, one must bend several times.

\subsection{Further bending steps}

We shall use compact hyperbolic manifolds which contain several
disjoint separating totally geodesic hypersurfaces. Again, such
manifolds exist in all dimension, see \cite{Mill}. In low
dimensions, a vast majority of known examples of compact hyperbolic
manifolds have this property (they fall into infinitely many
distinct commensurability classes, see \cite{Allcock}). Given such a
manifold, bending can be performed several times in a row. The next
lemmas show that at each step, the Zariski closure strictly
increases.

\begin{lemma}
\label{induction} Let $X'=H_{\bh}^{n}$. Let $Z$ be a totally
geodesic subspace of $X=H_{\bh}^{m}$ such that $X'\subsetneq
Z\subsetneq X$. Then $Z$ is a quaternionic subspace. Furthermore,
there exists an $a\in G$ fixing $X'$ pointwise which does not map
$Z$ into itself.
\end{lemma}

\begin{pf}
Otherwise, $Z$ would be $Sp(m-n)$-invariant. In particular, for
$x\in X'$, $T_x Z$ would be $Sp(m-n)$-invariant. Since $Sp(m-n)$
acts irreducibly on $(T_x X')^{\bot}$, $Z$ must be equal to $X'$ or
$X$, a contradiction. $Z$ is a negatively curved symmetric space
containing $H_{\bh}^n$, $n\geq 2$, so it is a quaternionic subspace.
\end{pf}

\begin{Prop}
\label{bent} Let $M$ be a compact hyperbolic $n$-manifold. Let
$m\geq n$. Assume that $M$ contains $N$ disjoint separating totally
geodesic hypersurfaces. Let $\Gamma=\pi_1 (M)\subset PO(n,1)^0 \ra
PSp(m,1)$. If $N\geq m-n+2$, then $\Gamma$ can be continuously
deformed to a Zariski dense subgroup of $PSp(m,1)$.
\end{Prop}

\begin{pf}
According to Corollary \ref{complex}, a first bending in $PU(n,1)$
provides us with a Zariski dense subgroup of $PU(n,1)$.

A second bending in $PSp(n,1)$ gives a Zariski dense subgroup of
$PSp(n,1)$. Indeed, the fixator of $H_{\br}^{n-1}$ is an
$Sp(1)Sp(1)$ which contains an element $a$ which does not map
$H_{\bc}^{n}$ to itself. By Lemma \ref{totallycomplex}, no proper
totally geodesic subspace of $H_{\bh}^{n}$ contains both
$H_{\bc}^{n}$ and $a(H_{\bc}^{n})$. Lemma \ref{zariski} implies that
the bent subgroup is Zariski dense.

A third series of bendings gives a Zariski dense subgroup of
$PSp(m,1)$. Lemma \ref{induction} allows inductively to select a
parameter $a$ which strictly increases the dimension of the totally
geodesic hull. After at most $m-n$ more steps, the obtained subgroup
is Zariski dense, thanks to Lemma \ref{zariski}.
\end{pf}

\subsection{Bending along laminations}

Since we need to bend surfaces of genus as low as 2, which do not
admit pairs of disjoint separating closed geodesics, we describe W.
Thurston's general construction of bending along totally geodesic
laminations, which does not require the leaves to be separating. We
stick to the special case of totally real, totally geodesic 2-planes
of $H_{\bh}^2$.

Let $Y=H_{\br}^{2}\subset H_{\bh}^2 =X$. If $\ell\subset Y$ is a
geodesic, the subgroup $Fix(\ell)$ of $Isom(X)$ that fixes $\ell$
pointwise is conjugate to $Sp(1)Sp(1)$. The Lie algebras of these
subgroups form an $\I\bh\oplus\I\bh$-bundle $\mathcal{B}$ over the
space $\mathcal{L}$ of geodesics in $Y$. Pick once et for all an
arbitrary Borel trivialization of this bundle. A {\em lamination} on
$Y$ is a closed subset of $\mathcal{L}$ consisting of pairwise non
intersecting geodesics. A {\em measured lamination} on $Y$ is the
data of a lamination $\lambda$ and a transverse
$\I\bh\oplus\I\bh$-valued measure. By a transverse measure, we mean
the data, for each continuous curve $c:[a,b]\to Y$ which crosses all
geodesics of $\lambda$ in the same direction, of a finite Borel
$\I\bh\oplus\I\bh$-valued measure $\mu_c$ on $[a,b]$, with the
following compatibility : if a curve $c':[a,b]\to Y$ can be deformed
to $c$ by sliding along $\lambda$, then $\mu_{c'}=\mu_c$. A discrete
collection of geodesics, with an $\I\bh\oplus\I\bh$-valued Dirac
mass at each geodesic, is a simple example of a measured lamination.
Since only such laminations will ultimately be used, we shall not
discuss non discrete measured laminations further.

The Lie algebra bundle $\mathcal{B}$ is a subbundle of the trivial
bundle with fiber the Lie algebra $\mathfrak{sp}(2,1)$. Therefore,
for every transversal curve $c$, the measure $\mu_c$ can be pushed
forward to yield an $\mathfrak{sp}(2,1)$-valued measure on $[a,b]$.
This measure integrates into a continuous map $[a,b]\to Sp(2,1)$,
see for example \cite{Es}. We denote the resulting element of
$Sp(2,1)$ by $\int \mu_c$. If $c=c_1 c_2$ is obtained by traversing
a first curve $c_1$ and then a second curve $c_2$, then Chasles rule
$\int\mu_{c_1 c_2}=(\int\mu_{c_1})(\int\mu_{c_2})$ holds, which
allows to extend the definition to curves which are piecewise
transversal. Define a map $f:Y\ra X$ as follows. Pick an origin
$o\in Y$. Given $y\in Y$, join $o$ to $y$ with a piecewise
transversal curve $c_y$ and set $f(y)=(\int \mu_{c_y})y$. One checks
that $f(y)$ does not depend on the choice of piecewise transversal
curve.

For instance, in the case of a discrete lamination, $f$ is piecewise
isometric and totally geodesic away from the support of $\lambda$.
At each geodesic $\ell$ of the lamination, $f$ bends, i.e. the
totally geodesic pieces of the surface $f(Y)$ at either side of
$\ell$ meet at a $Fix(\ell)$-angle equal to $\exp(\mu(\ell))$. The
general case is best understood by considering limits of discrete
measured laminations.

Let $\rho:\Gamma\ra Sp(2,1)$ be an isometric action of a group
$\Gamma$ which leaves $Y$ and the measured lamination invariant.
Then, for every piecewise transversal curve $c$, and
$\gamma\in\Gamma$,
$\int\mu_{\rho(\gamma)(c)}=\rho(\gamma)(\int\mu_c)\rho(\gamma)^{-1}$.
For $\gamma\in\Gamma$, let $\rho_{\lambda}(\gamma)=(\int
\mu_{c_{\gamma}})\rho(\gamma)$, where $c_{\gamma}$ is a piecewise
transversal curve joining $o$ to $\rho(\gamma)o$. Then
$\rho_{\lambda}:\Gamma\ra Sp(2,1)$ is a homomorphism which
stabilizes $f(Y)$, and $f$ is equivariant. Indeed, let $c_1$ (resp.
$c_2$) be a piecewise transversal curve joining $o$ to
$\rho(\gamma_1 )o$ (resp. to $\rho(\gamma_2 )o$). Then $c_1
\rho(\gamma_1 )(c_2)$ joins $o$ to $\rho(\gamma_1 \gamma_2 )o$ and
\begin{eqnarray*}
\rho_{\lambda}(\gamma_1 \gamma_2)&=&(\int\mu_{c_1 \rho(\gamma_1 )(c_2)})\rho(\gamma_1 \gamma_2) \\
&=&(\int\mu_{c_1})(\int\mu_{\rho(\gamma_1 )(c_2)})\rho(\gamma_1 \gamma_2 )\\
&=&(\int\mu_{c_1})\rho(\gamma_1 )(\int\mu_{c_2})\rho(\gamma_1^{-1}) \rho(\gamma_1 \gamma_2 )\\
&=&\rho_{\lambda}(\gamma_1) \rho_{\lambda}(\gamma_2).
\end{eqnarray*}
If $y\in Y$ and $\gamma\in\Gamma$, let $c_{y}$ (resp. $c_{\gamma}$)
be a piecewise transversal curve joining $o$ to $y$ (resp. to
$\rho(\gamma)o$). Then $c_{\gamma}\rho(\gamma)(c_{y})$ joins $o$ to
$\rho(\gamma)y$, thus
\begin{eqnarray*}
f(\rho(\gamma)y)&=&(\int\mu_{c_{\gamma}\rho(\gamma)(c_{y})})\rho(\gamma)y\\
&=&(\int\mu_{c_{\gamma}})(\int\mu_{\rho(\gamma)(c_{y})})\rho(\gamma)y\\
&=&(\int\mu_{c_{\gamma}})\rho(\gamma)(\int\mu_{c_{y}})\rho(\gamma)^{-1}\rho(\gamma)y\\
&=&\rho_{\lambda}(\gamma)f(y).
\end{eqnarray*}

\begin{Prop}
\label{lamination} Let $\Sigma$ be a closed hyperbolic surface with
fundamental group $\Gamma$. Map $\Gamma\ra SO(2,1)\ra Sp(2,1)$.
There exist measured laminations $\lambda$ on $\Sigma$ which make
the bent group $\rho_{\lambda}(\Gamma)$ Zariski dense in $Sp(2,1)$.
\end{Prop}

\begin{pf}
As a lamination, take the lifts to $Y=\tilde{\Sigma}$ of two
disjoint closed geodesics in $\Sigma$. A transversal measure in this
case is simply the data of elements $a_j \in Fix(\ell_j)$ for two
lifts $\ell_1$, $\ell_2$. Note that the components of the complement
of the two geodesics in $\Sigma$ are not simply connected. In other
words, each component of the complement of the support of the lifted
lamination on $Y$ is stabilized by a subgroup of $\Gamma$ which is
Zariski dense in $SO(2,1)$. It follows that the Zariski closure of
$\rho_{\lambda}(\Gamma)$ contains $SO(2,1)$. It also contains the
conjugates of $SO(2,1)$ by the two isometries $a_1$ and $a_2$.

According to Lemma \ref{zariski}, the Zariski closure of
$\rho_{\lambda}(\Gamma)$ contains the isometry group of the totally
geodesic hull $Z$ of $Y\cup a_1 (Y)\cup a_2 (Y)$. As in the proof of
Proposition \ref{bent}, bending by $a_1$ gives a group which is
Zariski dense in a conjugate of $PU(2,1)$, bending by $a_1$ and
$a_2$ gives a group which is Zariski dense in $PSp(2,1)$.
\end{pf}

\section{Flexibility of Fuchsian surface groups}
\label{fuchsian}

In this section, we investigate homomorphisms of a surface group
into $Sp(2,1)$ in a neighborhood of the embedding via $SU(1,1)$ and
$Sp(1,1)$. We shall call them {\em Fuchsian}, to distinguish them
from the bendable homomorphisms arising from the embedding via
$SO(2,1)$.

\subsection{Second order calculations}
Let $S$ be a compact Riemann surface with  genus $>1$ and
$\rho_0:\pi_1(S)=\Gamma \subset SU(1,1)\ra  Sp(1,1)\subset Sp(2,1)$
be a standard representation fixing a quaternionic line in
$H^2_\bh$. Since $H^1(\pi_1 (S),\bh^2)\not=0$, Proposition
\ref{relative} does not apply. We have to investigate which
infinitesimal deformations represented by $H^1(\pi_1
(S),\mathfrak{sp}(2,1))$ are integrable.

The second order integrability condition for infinitesimal
deformations at $\phi$ of representations of a group $\Gamma$ in a
Lie group $G$ can be expressed in terms of the {\em cup-product}, a
symmetric bilinear map
\begin{eqnarray*}
[\cdot ,\cdot]:H^1(\Gamma, \mathfrak{g}_{Ad\phi})\to H^2(\Gamma,
\mathfrak{g}_{Ad\phi}).
\end{eqnarray*}
For $u\in Z^1(\Gamma, \mathfrak{g}_{Ad\phi})$,
$$
[u,u](\alpha,\beta)=[u(\alpha),Ad\phi(\alpha)u(\beta)].
$$
It is well-known, \cite{NR}, that for a representation $\phi$ from
$\Gamma$ to a reductive group $G$, if there exists a smooth path
$\phi_t$ in $Hom(\Gamma, G)$ which is tangent to $u\in Z^1(\Gamma,
\mathfrak{g}_{Ad\phi})$, then $[u,u]=0$. According to Theorem 3 in
\cite{G1}, for surface groups, this necessary condition is also
sufficient.

\begin{thm}
\label{Goldman} {\em (W. Goldman)}. Let $S$ be a closed surface, let
$G$ be a reductive group. Let $\phi:\pi_1(S) \ra G$ be a
representation such that the Zariski closure of $\phi(\pi_1(S))$ is
also reductive. Then for any $u\in
Z^1(\pi_1(S),\mathfrak{g}_{Ad\phi})$, $[u,u]=0$ if and only if there
exists an analytic path $t\mapsto\phi_t$ in $Hom(\pi_1(S),G)$ which
is tangent to $u$.
\end{thm}

\subsection{Splitting of the cup-product map}

The centralizer of $SU(1,1)$ in $Sp(2,1)$ is $Sp(1)\times U(1)$,
where $Sp(1)$ is the centralizer of $Sp(1,1)$ and $U(1)\subset
Sp(1,1)$ is the centralizer of $SU(1,1)$ in $Sp(1,1)$. Then by
Poincar\'e duality
$$
H^2 (\pi_1(S),\mathfrak{sp}(2,1))=H^0 (\pi_1(S),\mathfrak{sp}(2,1))=
\mathfrak{sp}(1)\oplus \mathfrak{u}(1).
$$

Let $u\in H^1(\pi_1(S),\mathfrak{sp}(2,1))$ split as $u=u_{\mathfrak{sp}(1)}+u_{\mathfrak{sp}(1,1)}+u_{\bh^2}$. Since $\mathfrak{sp}(1,1)$ and $\mathfrak{sp}(1)$ commute, $[u_{\mathfrak{sp}(1)},u_{\mathfrak{sp}(1,1)}]=0$. Since the subspace $\bh^2 \subset \mathfrak{sp}(2,1)$ is $Sp(1)\times Sp(1,1)$-invariant, $[u_{\mathfrak{sp}(1)},u_{\bh^2}]$ and $[u_{\mathfrak{sp}(1,1)},u_{\bh^2}]$ belong to $H^2(\pi_1(S),\bh^2)=0$. Therefore
$$
[u,u]=[u_{\mathfrak{sp}(1)},u_{\mathfrak{sp}(1)}]+[u_{\mathfrak{sp}(1,1)},u_{\mathfrak{sp}(1,1)}]+[u_{\bh^2},u_{\bh^2}].
$$
Since $\mathfrak{sp}(1)$ and $\mathfrak{sp}(1,1)$ are subalgebras, $[u_{\mathfrak{sp}(1)},u_{\mathfrak{sp}(1)}]$ belongs to $H^2(\pi_1(S),\mathfrak{sp}(1))=\mathfrak{sp}(1)$, and $[u_{\mathfrak{sp}(1,1)},u_{\mathfrak{sp}(1,1)}]$ belongs to $H^2(\pi_1(S),\mathfrak{sp}(1,1))=\mathfrak{u}(1)$. On the other hand, $[u_{\bh^2},u_{\bh^2}]$ has nontrivial components $[u_{\bh^2},u_{\bh^2}]_{\mathfrak{u}(1)}$ and $[u_{\bh^2},u_{\bh^2}]_{\mathfrak{sp}(1)}$ on both $H^2(\pi_1(S),\mathfrak{sp}(1,1))$ and $H^2(\pi_1(S),\mathfrak{sp}(1))$.

\subsection{Homomorphisms to $Sp(1)$}

In the special case of the trivial representation to $Sp(1)$, the cup-product map can be computed.

\begin{lemma}
\label{sp1} Let $S$ be a closed surface. Let $\pi_1
(S)$ act trivially on $\mathfrak{sp}(1)$. The quadratic map $H^1 (\pi_1 (S),\mathfrak{sp}(1))\to H^2 (\pi_1 (S),\mathfrak{sp}(1))$, $u\mapsto [u,u]$, is onto.
\end{lemma}

\begin{pf}
Here, $H^1(\pi_1(S),\mathfrak{sp}(1))\simeq
H^1(\pi_1(S),\br)\otimes\mathfrak{sp}(1))$, If $a$, $b\in
H^1(\pi_1(S),\br)$ and $q$, $q'\in \mathfrak{sp}(1)$, then
\begin{eqnarray*}
[a\otimes q,b\otimes q']=a\smile b\otimes [q,q'].
\end{eqnarray*}
For every $q''\in \mathfrak{sp}(1)$, there exist $q$, $q'\in \mathfrak{sp}(1)$ such that $[q,q']=q''$. Poincar\'e duality implies that there exist $a$, $b\in
H^1(\pi_1(S),\br)$ such that $a\smile b\not=0$. Therefore the cup-product map is onto. 
\end{pf}

\subsection{Homomorphisms to $Sp(1,1)$}

A similar statement applies to $H^1(\pi_1(S),\mathfrak{sp}(1,1))$.

\begin{lemma}
\label{sp11} Let $S$ be a closed hyperbolic surface. View $\pi_1
(S)$ as a subgroup of $SU(1,1)\subset Sp(1,1)$. The quadratic map $H^1 (\pi_1 (S),\mathfrak{sp}(1,1))\to H^2 (\pi_1 (S),\mathfrak{sp}(1,1))=\mathfrak{u}(1)$, $u\mapsto [u,u]$, is onto.
\end{lemma}

\begin{pf}
$\mathfrak{sp}(1,1)$ consists of quaternionic $2\times 2$ matrices $\begin{pmatrix}
a & b  \\
\bar{b} & d
\end{pmatrix}$ with $a$, $d$ imaginary quaternions. The complex matrices in $\mathfrak{sp}(1,1)$ form the subalgebra $\mathfrak{u}(1,1)=\mathfrak{su}(1,1)\oplus\mathfrak{u}(1)$, where $\mathfrak{u}(1)$ consists of complex imaginary multiples of the unit matrix. As a $U(1,1)$-invariant projection $\mathfrak{sp}(1,1)\to\mathfrak{u}(1)=\br$, we can use the linear form
\begin{eqnarray*}
\pi_{\mathfrak{u}(1)}\begin{pmatrix}
a & b  \\
\bar{b} & d
\end{pmatrix}=\Re e(i(a+d)).
\end{eqnarray*}
Let $W$ denote the set of matrices of the form $j\begin{pmatrix}
z & w  \\
-w & t
\end{pmatrix}$, where $z$, $w$ and $t\in\bc$. Then $W$ is a $U(1,1)$-invariant complement of $\mathfrak{u}(1,1)$ in $\mathfrak{sp}(1,1)$. Given two elements $X=j\begin{pmatrix}
z & w  \\
-w & t
\end{pmatrix}$ and $X'=j\begin{pmatrix}
z' & w'  \\
-w' & t'
\end{pmatrix}$ in $W$, one computes
\begin{eqnarray*}
\pi_{\mathfrak{u}(1)}([X,X'])=2\I(\bar{z}z'+\bar{t}t'-2\bar{w}w').
\end{eqnarray*}
This is a symplectic structure on $W$ (viewed as a real vector space). From Poincar\'e duality for local coefficient systems, it follows that the quadratic form $\pi_{\mathfrak{u}(1)}([\cdot,\cdot])$ on $H^1 (\pi_1 (S),W)$ is nondegenerate. In particular, it is onto. {\em A fortiori}, the quadratic form $[\cdot,\cdot]$ on $H^1 (\pi_1 (S),\mathfrak{sp}(1,1))$ is onto.
\end{pf}

\subsection{Flexibility of certain Fuchsian surface groups}

A surface group in $SU(n,1)$ is {\em Fuchsian} if it  stabilizes a
complex line in complex hyperbolic space. Let us extend the
terminology. Say a surface group in $Sp(n,1)$ is {\em Fuchsian} if
it stabilizes a complex line in quaternionic hyperbolic space. Note
that every complex line is contained in a unique quaternionic line.

It is well-known that Fuchsian groups in $SU(2,1)$ (or, more
generally, $SU(n,1)$) cannot be deformed to Zariski dense groups. We
show that when $SU(2,1)$ is embedded in the larger group $Sp(2,1)$,
this rigidity property fails. We make essential use of the main
result of \cite{G1}.

\begin{Prop}
\label{surface1} Let $S$ be a compact Riemann surface with  genus
$>1$ and $\rho_0:\pi_1(S)=\Gamma \subset SU(1,1)\ra  Sp(1,1)\subset
Sp(2,1)$ be a standard representation fixing a quaternionic line in
$H^2_\bh$. Then there exist local deformations of $\rho_0$ which do
not stabilize any quaternionic line.
\end{Prop}

\begin{pf}
Let $u\in H^1(\pi_1 (S),\bh^2)$ be nonzero. According to Lemmas
\ref{sp1} and \ref{sp11}, there exist $v\in H^1(\pi_1 (S),\mathfrak{sp}(1))$ and $w\in H^1(\pi_1 (S),\mathfrak{sp}(1,1))$ such that $[v,v]=-[u,u]_{\mathfrak{sp}(1)}$ and $[w,w]=-[u,u]_{\mathfrak{u}(1)}$. Then $x=u+v+w\in H^1(\pi_1 (S),\mathfrak{g})$ is nonzero and satisfies $[x,x]=0$. According to Goldman's Theorem
\ref{Goldman}, there exists an analytic curve $t\mapsto\rho_t$ in
$Hom(\pi_1 (S),G)$, starting at $\rho_0$, whose initial speed is a
representative of the cohomology class $x$. Since $x\notin H^1(\pi_1
(S),\mathfrak{sp}(1)\oplus \mathfrak{sp}(1,1))$, for $t\not=0$
small, $\rho_t$ cannot be conjugated to the subgroup $Sp(1,1)Sp(1)$,
i.e., does not stabilize any quaternionic line.
\end{pf}

\noindent\textbf{Proof of Theorem \ref{2}}. Proposition
\ref{surface1} is statement (2) of Theorem \ref{2}. Statement (1) of
Theorem \ref{2} is a consequence of the bending construction. For
surfaces of sufficiently high genus, one can apply Proposition
\ref{bent}. In low genus, one needs bend along a geodesic
lamination, see Proposition \ref{lamination}.

\section{Discrete representations}
\label{discrete}

\begin{Prop}
Let $\Gamma$ be a uniform lattice in $Sp(1,1)$. Let $\rho:\Gamma\ra
Sp(2,1)$ be a discrete and faithful homomorphism. Then,
\begin{itemize}
  \item either $\rho$ is standard, i.e. it stabilizes a quaternionic line,
  \item or the image is Zariski dense.
\end{itemize}
\end{Prop}

\begin{pf}
Suppose $\rho(\Gamma)$ is not Zariski dense. Then it cannot be
contained in a parabolic subgroup of $Sp(2,1)$ since $\Gamma$ is not
solvable. So it must stabilizes a totally geodesic subspace of
$H^2_\bh$, see \cite{Kim1}. If it stabilizes a quaternionic line, it
is a standard representation, by Mostow rigidity. Suppose it
stabilizes $H^2_\bc$. Then $H^2_\bc/\rho(\Gamma)$ is a manifold. If
it is not closed, the cohomological dimension of $\Gamma$ cannot be
4, which contradicts $\Gamma$ being a uniform lattice in $Sp(1,1)$.
So $H^2_\bc/\rho(\Gamma)$ is a closed manifold, which implies that
$H^2_\bc$ and $H^4_\br$ are quasi-isometric, which is impossible,
again by a result of G.D. Mostow.
\end{pf}

We suspect that there is no Zariski dense discrete faithful group
$\rho(\Gamma)$. \\

\bigskip

{\bf Acknowledgement}. We thank an anonymous referee for his valuable
suggestions.

1991 {\sl{Mathematics Subject Classification.}}51M10, 57S25.

{\sl{Key words and phrases.}} Quaternionic hyperbolic space, rank
one symmetric space, quasifuchsian representation, bending,
rigidity, group cohomology

\vskip1cm

\noindent     Inkang Kim\\
     Department of Mathematics\\
     Seoul National University\\
     Seoul, 151-742, Korea\\
     \texttt{inkang\char`\@ math.snu.ac.kr}

\smallskip

     \noindent  Pierre Pansu\\ Laboratoire de Math{\'e}matiques
     d'Orsay\\
     UMR 8628 du CNRS\\
 Universit{\'e} Paris-Sud\\
 91405 Orsay C\'edex, France\\
  \texttt{pierre.pansu\char`\@ math.u-psud.fr}

     \end{document}